\documentclass[11pt]{article}
\usepackage{amsmath, epsfig, cite,colordvi,xcolor}
\usepackage{graphicx,subfigure}
\usepackage{amssymb}
\usepackage{amsfonts}
\usepackage{multirow}
\usepackage{array}
\usepackage{latexsym}
\usepackage{graphicx}
\usepackage{epstopdf}
\usepackage{tikz}
\usepackage{float}
\usepackage{cases}
\usepackage{mathrsfs}[12pt]
\usepackage{bm, amscd,  amsfonts, amssymb, cite,texdraw}
\usepackage{amsmath, epsfig, cite, color, colordvi}
\usepackage{amsmath,amsfonts,amssymb,amsthm}
\usepackage[colorlinks,urlcolor=purple,linkcolor=red,anchorcolor=green,citecolor=blue]{hyperref}
\usepackage{enumerate}
\usepackage{titlesec}
\usepackage{appendix}
\theoremstyle{plain}
\newtheorem{thm}{Theorem}[section]
\newtheorem{lem}[thm]{Lemma}

\newtheorem{conj}[thm]{Conjecture}

\numberwithin{equation}{section}

\def\lnk{\mathcal{L}_{n,k}^{(r)}}

\def\qed{\hfill \rule{4pt}{7pt}}
\def\pf{\noindent {\it{Proof.} \hskip 2pt}}

\def\red{}

\begin{document}

\begin{center}

{\bf\Large The Tur\'{a}n problem for a family of tight linear forests}\\[20pt]

{\large  Jian Wang\footnote{\ Corresponding author:
wangjian01@tyut.edu.cn}, Weihua Yang} \\[10pt]

Department of Mathematics,\\
Taiyuan University of Technology, Taiyuan 030024, P.R. China.
\end{center}

\noindent {\bf Abstract.}
Let $\mathcal{F}$ be a family of $r$-graphs. The Tur\'{a}n number $ex_r(n;\mathcal{F})$
is defined to be the maximum number of edges in an $r$-graph of order $n$ that is $\mathcal{F}$-free. The famous Erd\H{o}s Matching Conjecture shows that
\[
ex_r(n,M_{k+1}^{(r)})= \max\left\{\binom{rk+r-1}{r},\binom{n}{r}-\binom{n-k}{r}\right\},
\]
where $M_{k+1}^{(r)}$ represents the $r$-graph consisting of $k+1$ disjoint edges. Motivated by this conjecture, we consider the Tur\'{a}n problem for tight linear forests. A tight linear forest is an $r$-graph whose connected components are all tight paths or isolated vertices.
Let $\mathcal{L}_{n,k}^{(r)}$ be the family of all tight linear forests of order $n$ with $k$ edges in $r$-graphs. In this paper, we prove that for sufficiently large $n$,
\[
ex_r(n;\mathcal{L}_{n,k}^{(r)})=\max\left\{\binom{k}{r},
\binom{n}{r}-\binom{n-\left\lfloor (k-1)/r\right \rfloor}{r}\right\}+d,
\]
where $d=o(n^r)$ and if $r=3$ and $k=cn$ with $0<c<1$; if $r\geq 4$ and $k=cn$ with $0<c<1/2$. The proof is based on the  weak regularity lemma for hypergraphs. We also conjecture that for arbitrary $k$ satisfying $k \equiv 1\ (mod\ r)$, the error term $d$ in the above result equals 0. We prove that the proposed conjecture implies the Erd\H{o}s Matching Conjecture directly.


\allowdisplaybreaks

\section{Introduction}

Given $r\geq 2$, an $r$-graph  (or $r$-uniform hypergraph) $H$ is a pair $H=(V,E)$, where $V$ is a finite vertex set and $E$ is a family of $r$-element subsets of $V$. Let $\mathcal{F}$ be a family of $r$-graphs. A hypergraph $H$ is called
$\mathcal{F}$-free if for any $F\in \mathcal{F}$, there is no subgraph of $H$ isomorphic to $F$. The {\it Tur\'{a}n number}
$ex_r(n;\mathcal{F})$ is defined to be the maximum number of edges in an $r$-graph on $n$ vertices
that is $\mathcal{F}$-free. For a single $r$-graph $F$, we write $ex_r(n; F)$ instead of $ex_r(n; \{F\})$. Tur\'{a}n
introduced this problem in \cite{turan1}, and we recommend \cite{sidorenko,keevash} for surveys on Tur\'{a}n problems for graphs and hypergraphs.

Let $H$ be an $r$-graph on $n$ vertices. A {\it matching} $M$ in $H$ is a collection of disjoint edges of $H$. We denote by $\nu(H)$ the number of edges in a maximum matching of $H$. The following classical conjecture
is due to Erd\H{o}s \cite{erdos4} proposed in 1965.
\begin{conj}\label{matchingconj}[Erd\H{o}s Matching Conjecture, \cite{erdos4}]
Let $H$ be an $r$-graph on $n$ vertices with $\nu(H)=k$. Then
\[
e(H)\leq \max\left\{\binom{rk+r-1}{r},\binom{n}{r}-\binom{n-k}{r}\right\}.
\]
\end{conj}

Let $H_1$ and $H_2$ be two disjoint $r$-graphs. The {\it join} of two $r$-graphs, denoted by $H_1\vee H_2$,
is defined as $V(H_1\vee H_2)=V(H_1)\cup V(H_2)$ and
$E(H_1\vee H_2)=E(H_1)\cup E(H_2)\cup \{S\in \binom{V(H_1)\cup
V(H_2)}{r}\colon S\cap V(H_1)\neq \emptyset \mbox{ and } S\cap V(H_2)\neq\emptyset\}$.
We denote by $K_n^{(r)}$ and $E_n^{(r)}$ the complete $r$-graph of order $n$ and the empty $r$-graph of order $n$, respectively. The constructions $K_{rk+r-1}^{(r)}$ and $K_{k}^{(r)}\vee E_{n-k}^{(r)}$ implies that the bounds  given in the Conjecture \ref{matchingconj}, if true, is tight.

Let $M_{k}^{(r)}$ be an $r$-graph with exact $k$ disjoint edges. The conclusion of the Erd\H{o}s Matching Conjecture can also be expressed as a Tur\'{a}n function as follows:
\[
ex_r(n,M_{k+1}^{(r)})= \max\left\{\binom{rk+r-1}{r},\binom{n}{r}-\binom{n-k}{r}\right\}.
\]
The case $k=1$ for Erd\H{o}s Matching Conjecture is the classical Erd\H{o}s-Ko-Rado Theorem\cite{erdos3}.
For $r=1$ the conjecture is trivial and for $r=2$ it was proved by Erd\H{o}s and Gallai\cite{erdos2}.
When $r=3$, Frankl, R\"{o}dl and Ruci\'{n}ski\cite{frankl0} proved the conjecture for $k\leq n/4$.
Recently, Frankl\cite{frankl2} proved the conjecture for $r=3$. For arbitrary $r$, Bollob\'{a}s, Daykin and Erd\H{o}s\cite{bollobas2} proved the conjecture
 for $n>2r^3(k-1)$.
Huang, Loh and Sudakov\cite{huang} improved it to $n\geq 3r^2k$. Recently, Frankl\cite{frankl1} proved
the conjecture for $n\geq (2r-1)k+r$.

A matching can also be viewed  as a forest of paths with length one. In general, we can consider forests of paths with given
length. For graphs, let $P_t$ denote a path on $t$ vertices, and $k\cdot P_t$ denote $k$ vertex-disjoint copies of $P_t$.
Bushaw and Kettle\cite{bushaw} determined the exact values of $ex(n,k\cdot P_t)$ for $n$ appropriately large relative to $k$ and $t$. Later, Yuan and Zhang\cite{zhang} determined the value $ex(n, k \cdot P_3)$ and characterize all extremal
graphs for all $k$ and $n$.
Denote by
$L(s_1,\ldots,s_k)$ a linear forest consisting of $k$ vertex-disjoint paths with $s_1, \ldots, s_k$ edges. Suppose every $s_i\geq 1$ and at least one $s_i\neq 2$, then
Lidicky, Liu, and Palmer\cite{lidicky} proved that for sufficiently large $n$,
\begin{align*}
ex(n;L(s_1,\ldots,s_k)) = \binom{n}{2}-\binom{n-t}{2}+c',
\end{align*}
where $t=\sum_{i=1}^k\left\lfloor\frac{s_i+1}{2}\right\rfloor -1$, and $c' = 1$ if all $s_i$ are even and $c' = 0$ otherwise.

Instead of one single linear forest, Wang and Yang\cite{wj} consider a family of linear forests as forbidden subgraphs. For the purpose of simplification, they allow isolated vertices in a linear forest in their paper, i.e. a {\it linear forest} is a graph consisting of vertex-disjoint paths or isolated vertices.  Let $\mathcal{L}_{n,k}$ be the set of all linear forests of order $n$
with $k$ edges. The problem of determining the Tur\'{a}n number $ex(n;\mathcal{L}_{n,k})$ was introduced in \cite{wj}. Recently, Ning and Wang\cite{ningbo} prove that
\[
ex(n;\mathcal{L}_{n,k})=
\max \left\{\binom{k}{2},\binom{n}{2}-\binom{n-\left\lfloor (k-1)/2\right \rfloor}{2} + c\right\},
\]
where $c=0$ if $k$ is odd and $c=1$ otherwise.

In this paper, we generalize this problem into hypergraphs. Define a {\it tight linear forest} to be an $r$-graph consisting of
vertex-disjoint tight paths or isolated vertices. The {\it maximum tight linear forest} in $H$ is a subgraph
of $H$ with maximum number of edges that is a tight linear forest. We denote by $l(H)$ this maximum number.  When $r=2$, it reduced to the linear forest in graphs.
Let $\lnk$ be the family of all tight linear forests of order $n$ with $k$ edges. In this paper, we consider the Tur\'{a}n number of $\lnk$ and prove the following theorem.
\begin{thm}\label{main2}
Let $\lnk$ be the family of all tight linear forests of order $n$ with $k$ edges.
For $r\geq 3$, $n$ sufficiently large, $0<c<1$ when $r=3$, and $0<c<1/2$ when $r\geq 4$,
\[
ex_r(n;\mathcal{L}_{n,cn}^{(r)})=\max\left\{\binom{cn}{r},
\binom{n}{r}-\binom{n-\frac{cn}{r}}{r}\right\}+o(n^r).
\]
\end{thm}
We also propose a conjecture as follows:
\begin{conj}\label{conj} Let $\lnk$ be the family of all tight linear forests of order $n$ with $k$ edges. For $k>r$ and $k\equiv 1(mod\ r)$,
\[
ex_r(n;\lnk)=\max\left\{\binom{k+r-2}{r},\binom{n}{r}-\binom{n-(k-1)/r}{r}\right\}.
\]
\end{conj}

For $r=2$, the conjecture is true according to Ning and Wang's result\cite{ningbo}.  For $r=3$, the conjecture is asymptotically true for $k=cn$ with $0<c<1$ and $n$ sufficiently large, according to Theorem \ref{main2}. For $r\geq 4$, the conjecture is asymptotically true for $k=cn$ with $0<c<\frac{1}{2}$ and $n$ sufficiently large, according to Theorem \ref{main2}.

It should be noticed that the Conjecture \ref{conj} implies the Erd\H{o}s Matching Conjecture.
Let $H$ be an $r$-graph on $n$ vertices with  $\nu(H)=k$. Since each tight path contains a matching
with at least one $r$-th of its edges, each tight linear forest contains a matching with
at least one $r$-th of its edges. It follows that $l(H)\leq rk$. Therefore, by Conjecture \ref{conj} we have
\begin{align*}
e(H) &\leq ex(n;\mathcal{L}_{n,rk+1}^{(r)})\\  &=\max\left\{\binom{rk+1+r-2}{r},\red{\binom{n}{r}-\binom{n- (rk+1-1)/r}{r}}\right\}\\
&=\max\left\{\binom{rk+r-1}{r},\binom{n}{r}-\binom{n-k}{r}\right\},
\end{align*}
which is exactly the result given in the Erd\H{o}s Matching Conjecture.

The rest of the paper is organized as follows. In Section 2, we give two useful lemmas.
In Section 3, we introduce the weak regularity lemma for hypergraphs. In Section 4, we prove the Theorem \ref{main2} for $r\geq 4$. In Section 5, we prove the Theorem \ref{main2} for $r=3$. In Section 6, we give some concluding remarks.

\section{The matching lemma}

Recently, Frankl\cite{frankl1,frankl2} prove the follow two theorems, which confirms the Erd\H{o}s Matching Conjecture for $r\geq 3$ and $r\geq 4$ with $n\geq (2r-1)k+r$.

\begin{thm}\label{franklthm1}
Let $H$ be an $3$-graph on $n$ vertices with $\nu(H)=k$ and $n\geq 3k+2$. Then
\[
e(H)\leq \max\left\{\binom{3k+2}{3},\binom{n}{3}-\binom{n-k}{3}\right\}.
\]
\end{thm}

\begin{thm}\label{franklthm2}
Let $H$ be an $r$-graph on $n$ vertices with $\nu(H)=k$ and $n\geq (2r-1)k+r$. Then
\[
e(H)\leq \binom{n}{r}-\binom{n-k}{r}.
\]
\end{thm}
According to Theorem \ref{franklthm1} and \ref{franklthm2}, we can obtain lower bounds on the size of maximum matchings in dense hypergraphs. We call them the Matching Lemma.

\begin{lem}\label{matchinglem-3}
Let $H$ be an $3$-graph on $n$ vertices with $\alpha n^3$ edges where $n$ is sufficiently large and $\beta_0=\frac{\sqrt{321}-3}{52}\thickapprox 0.2869$. For $0< \alpha<\frac{9}{2}\beta_0^3$, $\nu(H)\geq (1-\sqrt[3]{1-6\alpha})n-2$ and for $\frac{9}{2}\beta_0^3< \alpha<\frac{1}{6}$,  $\nu(H)\geq \frac{\sqrt[3]{6\alpha}}{3}n-1$.
\end{lem}
\pf Let $\nu(H)=k=\beta n$. By Theorem \ref{franklthm1}, we have
\[
e(H)\leq \max\left\{\binom{3k+2}{3}, \binom{n}{3}-\binom{n-k}{3}\right\}.
\]
Clearly,
\[
\binom{3k+2}{3}=\frac{9}{2}\beta^3 n^3+o(n^3),
\mbox{ and } \binom{n}{3}-\binom{n-k}{3}=\frac{1-(1-\beta)^3}{6}n^3+o(n^3).
\]
Let $f(x)=27x^3+(1-x)^3-1$. Calculation shows that $\beta_0$ is a zero of $f(x)$ and
$f(x)<0$ for $x\in (0,\beta_0)$; $f(x)>0$ for $x\in (\beta_0,\frac{1}{3})$.
It follows that for  $\beta \in (0,\beta_0)$ and $n$ sufficiently large,
\[
\binom{3k+2}{3}< \binom{n}{3}-\binom{n-k}{3}.
\]
And for $\beta > \beta_0$ and $n$ sufficiently large,
\[
\binom{3k+2}{3}> \binom{n}{3}-\binom{n-k}{3}.
\]

For $0<\alpha<\frac{9}{2}\beta_0^3$, suppose to the contrary that $k< (1-\sqrt[3]{1-6\alpha})n-2$.
Since $\alpha<\frac{9}{2}\beta_0^3=\frac{1-(1-\beta_0)^3}{6}$, we have
$k=\beta n <(1-\sqrt[3]{1-6\alpha})n-2<\beta_0 n-2$. It follows that $\beta<\beta_0$. Therefore, by Theorem \ref{franklthm1}, we have
\begin{align*}
e(H)&\leq \binom{n}{3}-\binom{n-k}{3}\\
&< \binom{n}{3}-\binom{\sqrt[3]{1-6\alpha}n+2}{3}\\
&< \frac{n^3}{6}-\frac{(1-6\alpha)n^3}{6}=\alpha n^3,
\end{align*}
a contradiction. Thus, we conclude that $\nu(H)\geq (1-\sqrt[3]{1-6\alpha})n-2$ for $0< \alpha<\frac{9}{2}\beta_0^3$.

For $\frac{9}{2}\beta_0^3< \alpha<\frac{1}{6}$, suppose to the contrary that $k< \frac{\sqrt[3]{6\alpha}}{3}n-1$. If $\beta > \beta_0$,
then by Theorem \ref{franklthm1} we have
\begin{align*}
e(H)&\leq \binom{3k+2}{3} <\binom{\sqrt[3]{6\alpha}n-1}{3}<\alpha n^3,
\end{align*}
which leads to a contradiction. If $\beta < \beta_0$, then by Theorem \ref{franklthm1} we have
\begin{align*}
e(H)&\leq \binom{n}{3}-\binom{n-k}{3}\\
&= \frac{1-(1-\beta)^3}{6}n^3+o(n^3)\\
&< \frac{1-(1-\beta_0)^3}{6}n^3+o(n^3)\\
&= \frac{9\beta_0^3n^3}{2}+o(n^3)< \alpha n^3,
\end{align*}
a contradiction. Therefore, we conclude that $\nu(H)=k\geq \frac{\sqrt[3]{6\alpha}}{3}n-1$ for $\frac{9}{2}\beta_0^3< \alpha<\frac{1}{6}$.
\qed

\begin{lem}\label{matchinglem}
Let $H$ be an $r$-graph on $n$ vertices with $\alpha n^r$ edges where $n$ is sufficiently large.  If $0< \alpha<\frac{1-\left(1-\frac{1}{2r}\right)^r}{r!}$, then $\nu(H)\geq (1-\sqrt[r]{1-r!\alpha})n-r+1$.
\end{lem}

\pf  Let $\nu(H)=k=\beta n$.  Suppose to the contrary that $k< (1-\sqrt[r]{1-r!\alpha})n-r+1$. Since $\alpha<\frac{1-\left(1-\frac{1}{2r}\right)^r}{r!}$, it follows that
$k< (1-\sqrt[r]{1-r!\alpha})n-r+1<\frac{n}{2r}$. Thus, by Theorem \ref{franklthm2}, we have
\begin{align*}
e(H)&\leq \binom{n}{r}-\binom{n-k}{r}\\
&< \binom{n}{r}-\binom{\sqrt[r]{1-r!\alpha}n+r-1}{r}\\
&< \frac{n^r}{r!}-\frac{(1-r!\alpha)n^r}{r!}\\
&=\alpha n^r,
\end{align*}
a contradiction. Thus, we conclude that if $0< \alpha<\frac{1-\left(1-\frac{1}{2r}\right)^r}{r!}$, then $\nu(H)\geq (1-\sqrt[r]{1-r!\alpha})n-r+1$.\qed

\section{The weak regularity lemma for hypergraphs and path embeddings}

Szemer\'{e}di's Regularity Lemma\cite{szem} has been proved to be an incredibly powerful and useful tool in graph
theory as well as in Ramsey theory, combinatorial number theory and other areas of mathematics and theoretical
computer science. The weak regularity lemma for hypergraphs\cite{chung} is a straightforward extension of Szemer\'{e}di's Regularity Lemma.

Given an $r$-partite $r$-graph $H(X_1, X_2, \ldots, X_r)$. Define
\[
d(X_1, X_2, \ldots, X_r)=\frac{e(X_1, X_2, \ldots, X_r)}{|X_1||X_2|\cdots|X_r|}
\]
to be the edge density among $X_1$,$X_2$,\ldots,$X_r$.
Then $H(X_1, X_2, \ldots, X_r)$ is called $(\varepsilon,d(X_1, X_2, \ldots, X_r))$-regular ($\varepsilon$-regular, in short) if for all $r$-tuples of subsets $A_i\subset X_i$ with $|A_i|\geq \varepsilon|X_i|$ for all $i\in [r]$, we have
\[
|d(A_1, A_2, \ldots, A_r)-d(X_1, X_2, \ldots, X_r)|\leq \varepsilon.
\]
Now we state the weak regularity lemma for hypergraphs as follows.
\begin{thm}
Given $\varepsilon$, there exists $T=T(\varepsilon)$ and $n_0=n_0(\varepsilon)$ so that for every $r$-graph $H$ on $n$ vertices with $n\geq n_0$, there exists a partition $V=V_0\cup V_1\cup \ldots \cup V_t$ such that

(i) $\frac{1}{\varepsilon}\leq t\leq T$,

(ii) $|V_1|= |V_2|= \cdots = |V_t|$ and $|V_0|\leq \varepsilon n$,

(iii) for all but at most $\varepsilon t^r$ of $r$ tuples $(V_{i_1},V_{i_2},\ldots,V_{i_r})$, the induced $r$-partite $r$-graph among $V_{i_1}$,  $V_{i_2}$, $\ldots$, $V_{i_r}$ is $\varepsilon$-regular.
\end{thm}
The following lemma was essentially proved in \cite{han10}.
\begin{lem}\label{longpath}
Suppose $H$ is an $r$-partite $r$-graph with partition classes $V_1, V_2,\ldots,V_r$, $|V_i|=m$
for all $i\in [r]$, and $E(H)\geq d m^r$. Then there exists a tight path in $H$ with at least $\frac{dm}{2}$ edges.
\end{lem}

\begin{lem}\label{pathembedding}
Let $H$ be an $r$-partite $r$-graph with partition classes $V_1, V_2,\ldots,V_r$, $|V_i|=m$
for all $i\in [r]$. Suppose $H(V_1, V_2,\ldots,V_r)$ is $(\varepsilon,d)$-regular.
Then there are at most $\frac{3r}{(d-\varepsilon)\varepsilon}$ vertex disjoint tight
paths that cover all but at most $r\varepsilon m$ vertices of $H$.
\end{lem}

\pf
We greedily find disjoint tight paths of $rt$ vertices by Lemma \ref{longpath},
for some integer $t$. Since  $e(V_1, V_2,\ldots,V_r)= d m^r$, by Lemma \ref{longpath}
 we can find a tight path $P_1$ of $rt$ vertices with at least $\frac{dm}{3}$ edges for some integer $t$, such that for each $i$, $|V(P_1)\cap V_i|=t$. Then let $V_i^{(1)}=V_i\setminus V(P_1)$ for all $i \in [r]$.
   Let $m_1=|V_i^{(1)}|$. If $m_1\geq \varepsilon m$, then $d(V_1^{(1)}, V_2^{(1)},\ldots,V_r^{(1)})\geq d-\varepsilon$.
    It follows that $e(V_1^{(1)}, V_2^{(1)},\ldots,V_r^{(1)})\geq (d-\varepsilon)m_1^r$.
    Then we can find a tight path $P_2$  with at least $\frac{(d-\varepsilon)m_1}{3}$ edges.
    We can do it iteratively until $|V_i^{(n)}|<\varepsilon m$ for some $n$.
    Then we left at most $r\varepsilon m$ vertices. Since each tight path has
    at least $\frac{(d-\varepsilon)\varepsilon m}{3}$ vertices.
    Thus, we obtain at most $\frac{3r}{(d-\varepsilon)\varepsilon}$ vertex-disjoint tight paths that cover all but at most $r\varepsilon m$ vertices of $H$.
\qed

\section{The asymptotic Tur\'{a}n number of $\lnk$ for $r\geq 4$}

\begin{lem}\label{main}
Let $H$ be an $r$-graph on $n$ vertices with $\alpha n^r$ edges. For sufficiently large $n$ and $0< \alpha<\frac{1-\left(1-\frac{1}{2r}\right)^r}{r!}$, $H$ contains a tight linear forest with at least $(r-r\sqrt[r]{1-r!\alpha}-o(1))n$ edges.
\end{lem}
\pf Let $\varepsilon$ be a positive real number such that
\[\varepsilon \ll \min \left\{\alpha,\frac{r-r\sqrt[r]{1-r!\alpha}}{r^2+2r+r!(1-\frac{1}{2r})^{1-r}}\right\}.\]
Apply the weak regularity lemma to $H$ with the parameter $\varepsilon$ in a routing way. Then we get an $\varepsilon$-regular partition into $t+1$ sets $V_0,V_1,\ldots,V_t$ where $\frac{1}{\varepsilon}\leq t\leq T(\varepsilon)$. Then remove the following edges from $H$:
\begin{itemize}
   \item[(1)] edges with one endpoints in $V_0$.
   \item[(2)] edges intersecting $V_i$ for more than two vertices for $i=1,2,\ldots,t$.
   \item[(3)] edges among $V_{i_1}$,  $V_{i_2}$, $\ldots$, $V_{i_r}$ where $(V_{i_1},V_{i_2},\ldots,V_{i_r})$ is not $\varepsilon$-regular.
   \item[(4)] edges among $V_{i_1}$,  $V_{i_2}$, $\ldots$, $V_{i_r}$ where  $d(V_{i_1},V_{i_2},\ldots,V_{i_r})<2\varepsilon$.
 \end{itemize}
The number of edges in (1) is at most $\varepsilon n^r$. The number of edges in (2) is at most $t \left(\frac{n}{t}\right)^2n^{r-2}=\frac{n^r}{t}\leq \varepsilon n^r$. The number of edges in (3) is at most $\varepsilon t^r\left(\frac{n}{t}\right)^r = \varepsilon n^r$.  The number of edges in (4) is at most $\frac{t^r}{r!} 2\varepsilon \left(\frac{n}{t}\right)^r <\varepsilon n^r$. In total, we removed at most $4\varepsilon n^r$ edges. We call the remainder graph $H'$, and clearly $e(H')\geq (\alpha-4\varepsilon)n^r$.

Now we define a reduced $r$-graph $\Gamma$ on vertices set $\{1,2, \ldots,t\}$. If $(V_{i_1},V_{i_2},\ldots,V_{i_r})$ is $\varepsilon$-regular with density $d(V_{i_1},V_{i_2},\ldots,V_{i_r})\geq 2\varepsilon$, then $(i_1,i_2,\ldots,i_r)$ is an edge in $\Gamma$. Since $e_{H'}(V_{i_1},V_{i_2},\ldots,V_{i_r})>0$ if and only if $(i_1,i_2,\ldots,i_r)$ is an edge in $\Gamma$ and there are at most $\left(\frac{n}{t}\right)^r$ edges among $V_{i_1}$,  $V_{i_2}$, $\ldots$, $V_{i_r}$. Therefore,
\[
e(\Gamma)\geq \frac{e(H')}{\left(\frac{n}{t}\right)^r}\geq (\alpha-4\varepsilon)t^r.
\]
By Lemma \ref{matchinglem}, since $0< \alpha-4\varepsilon<\frac{1-\left(1-\frac{1}{2r}\right)^r}{r!}$, then $\nu(\Gamma)\geq (1-\sqrt[r]{1-r!(\alpha-4\varepsilon)})t-r+1$.

Let $M$ be a maximum matching in $\Gamma$. Then for each edge $(i_1,i_2,\ldots,i_r)\in M$, we can find at most $3r\varepsilon^{-2}$ vertex disjoint tight paths that cover all but at most $r\varepsilon \frac{n}{t}$ vertices of $V_{i_1}\cup V_{i_2}\cup \ldots\cup V_{i_r}$ by Lemma \ref{pathembedding}. Since all these tight paths cover at least $r\frac{(1-\varepsilon)n}{t}-r\varepsilon \frac{n}{t}=(1-2\varepsilon)\frac{rn}{t}$ vertices. Then there are at least $(1-2\varepsilon)\frac{rn}{t}-3r^2\varepsilon^{-2}\geq (1-3\varepsilon)\frac{rn}{t}$ edges in total.

Finally, by putting all the paths corresponding to each edge in $M$ together, we obtain that a tight linear forest $F$ with lots of edges. Thus, the number of edges in $F$ is at least
\begin{align*}
&\quad\left(\left(1-\sqrt[r]{1-r!(\alpha-4\varepsilon)}\right)t-r+1\right)\cdot(1-3\varepsilon)\frac{rn}{t}\\
&= \left(1-\sqrt[r]{1-r!(\alpha-4\varepsilon)}\right)(1-3\varepsilon)rn-(r-1)(1-3\varepsilon)\frac{rn}{t} \\
&\geq r\left(1-\sqrt[r]{1-r!(\alpha-4\varepsilon)}\right)n-3\varepsilon rn -(r-1)(1-3\varepsilon)\varepsilon rn \\
&\geq r\left(1-\sqrt[r]{1-r!\alpha}-(r-1)!(1-r!\alpha)^{-\frac{r-1}{r}}4\varepsilon\right)n-(r^2+2r)\varepsilon n\\
&\geq \left(r-r\sqrt[r]{1-r!\alpha}-\left(r^2+2r+r!\left(1-\frac{1}{2r}\right)^{1-r}\right)\varepsilon\right)n.
\end{align*}
where the first inequality follows from $t\geq \frac{1}{\varepsilon}$, the second inequality follows from Lagrange's Mean Value Theorem, and the last inequality follows from $\alpha \leq \frac{1-\left(1-\frac{1}{2r}\right)^r}{r!}$.

Thus, the Lemma holds.
\qed

\begin{thm}\label{mainr}
Let $\lnk$ be the family of all tight linear forests of order $n$ with $k$ edges.
For $r\geq 3$, $0<c<\frac{1}{2}$ and $n$ sufficiently large,
\[
ex_r(n;\mathcal{L}_{n,cn}^{(r)})=
\binom{n}{r}-\binom{n-\frac{cn}{r}}{r}+o(n^r).
\]
\end{thm}
\pf Let $\gamma$ be a positive real number. By Lemma \ref{main}, any $r$-graph on $n$ vertices with at least
$\left(\frac{1-\left(1-\frac{c}{r}\right)^r}{r!}+\gamma \right)n^r$
edges contains a linear forest with at least $\left(r-\sqrt[r]{(r-c)^r-r^rr!\gamma}-o(1)\right)n$\\$\geq cn$ edges. Thus, $ex(n;\mathcal{L}_{n,cn}^{(r)})\leq \left(\frac{1-\left(1-\frac{c}{r}\right)^r}{r!} \right)n^r+o(n^r)=\binom{n}{r}-\binom{n-\frac{cn}{r}}{r}+o(n^r)$. On the other hand, the construction $K_{(cn-1)/r}^{(r)}\vee E_{(1-c/r)n+1/r}^{(r)}$ implies that $ex(n;\mathcal{L}_{n,cn}^{(r)})= \binom{n}{r}-\binom{n-\frac{cn}{r}}{r}+o(n^r)$ for $c< 1/2$. Thus, we complete the proof.
\qed

\section{The asymptotic Tur\'{a}n number of $\mathcal{L}_{n,k}^{(3)}$}

\begin{lem}\label{main-3}
Let $H$ be an $3$-graph on $n$ vertices with $\alpha n^3$ edges and  $\beta_0=\frac{\sqrt{321}-3}{52}\thickapprox 0.2869$. For sufficiently large $n$ and $0< \alpha\leq\frac{9}{2}\beta_0^3$, $H$ contains a tight linear forest with at least $3(1-\sqrt[3]{1-6\alpha}-o(1))n$ edges; for $\frac{9}{2}\beta_0^3< \alpha<\frac{1}{6}$, $H$ contains a tight linear forest with at least $(\sqrt[3]{6\alpha}-o(1))n$ edges.
\end{lem}

\pf Let $\varepsilon$ be a positive real number such that
\[\varepsilon \ll \min \left\{\alpha,\frac{1-\sqrt[3]{1-6\alpha}}{13}\right\}\]
for $0< \alpha\leq\frac{9}{2}\beta_0^3$ and
\[\varepsilon \ll \min \left\{\alpha-\frac{9}{2}\beta_0^3,\frac{\sqrt[3]{6\alpha}}{21}\right\}\]
for $\frac{9}{2}\beta_0^3< \alpha<\frac{1}{6}$.
Apply the weak regularity lemma to $H$ with parameter $\varepsilon$ in a routing way, we get an $\varepsilon$-regular partition into $t+1$ sets $V_0,V_1,\ldots,V_t$ where $\frac{1}{\varepsilon}\leq t\leq T(\varepsilon)$. As the same as in the proof of Theorem \ref{main}, we get a reduced $3$-graph $\Gamma$ on vertices set $\{1,2, \ldots,t\}$ with at least $(\alpha-4\varepsilon)t^3$ edges. By Lemma \ref{matchinglem-3}, if $\frac{9}{2}\beta_0^3< \alpha-4\varepsilon<\frac{1}{6}$, then $\nu(\Gamma)\geq \frac{\sqrt[4]{6(\alpha-4\varepsilon)}}{3}t-1$; if $0< \alpha-4\varepsilon<\frac{9}{2}\beta_0^3$, then $\nu(\Gamma)\geq (1-\sqrt[3]{1-6(\alpha-4\varepsilon)})t-2$. Let $M$ be a maximum matching in $\Gamma$. Then for each edge $(i_1,i_2,i_3)\in M$, we can find at most $9\varepsilon^{-2}$ vertex disjoint tight paths that cover all but at most $3\varepsilon \frac{n}{t}$ vertices of $V_{i_1}\cup V_{i_2} \cup V_{i_3}$ by Lemma \ref{pathembedding}. Since all these tight paths cover at least $3\frac{(1-\varepsilon)n}{t}-3\varepsilon \frac{n}{t}=(1-2\varepsilon)\frac{3n}{t}$ vertices. Then there are at least $(1-2\varepsilon)\frac{3n}{t}-3\cdot 3^2\varepsilon^{-2}\geq (1-3\varepsilon)\frac{3n}{t}$ edges in total.

Finally, by putting all the paths corresponding to each edge in $M$ together, we obtain a tight linear forest $F$ with lots of edges. If $\frac{9}{2}\beta_0^3< \alpha-4\varepsilon<\frac{1}{6}$, then the number of edges in $F$ is at least
\begin{align*}
&\quad \left(\frac{\sqrt[3]{6(\alpha-4\varepsilon)}}{3}t-1\right)
\cdot (1-3\varepsilon)\frac{3n}{t}\\
& = (1-3\varepsilon)\sqrt[3]{6(\alpha-4\varepsilon)}n -(1-3\varepsilon)\frac{3n}{t}\\
&\geq(1-3\varepsilon)\sqrt[3]{6(\alpha-4\varepsilon)}n -(1-3\varepsilon)3\varepsilon n\\
&\geq(1-3\varepsilon)\left(\sqrt[3]{6\alpha}
-8\varepsilon\left(6(\alpha-4\varepsilon)\right)^{-\frac{2}{3}}\right)n-3\varepsilon n\\
&\geq \left(\sqrt[3]{6\alpha}-\left(3 + 3\sqrt[3]{6\alpha}+8
\left(6(\alpha-4\varepsilon)\right)^{-\frac{2}{3}}\right)\varepsilon\right)n\\
&\geq \left(\sqrt[3]{6\alpha}-\left( \frac{8}{9}\beta_0^{-2}+6\right)\varepsilon\right)n\\
& > \left(\sqrt[3]{6\alpha}-21\varepsilon\right)n,
\end{align*}
where the first inequality follows from $t\geq \frac{1}{\varepsilon}$, the second inequality follows from Lagrange's Mean Value Theorem, the third inequality follows from $\alpha -4\varepsilon \geq \frac{9}{2}\beta_0^3$ and $6\alpha<1$,  and the last third inequality follows from $\beta_0> \frac{1}{4}$.

If $0< \alpha -4\varepsilon<\frac{9}{2}\beta_0^3$, then $\nu(\Gamma)\geq (1-\sqrt[3]{1-6(\alpha-4\varepsilon)})t-2$, then the number of edges in $F$ is at least
\begin{align*}
&\quad\left(\left(1-\sqrt[3]{1-6(\alpha-4\varepsilon)}\right)t-2\right)\cdot(1-3\varepsilon)\frac{3n}{t}\\
&= \left(1-\sqrt[3]{1-6(\alpha-4\varepsilon)}\right)3(1-3\varepsilon)n-6(1-3\varepsilon)\frac{n}{t} \\
&\geq 3\left(1-\sqrt[3]{1-6(\alpha-4\varepsilon)}\right)n-9\varepsilon n -6(1-3\varepsilon)\varepsilon n \\
&\geq 3\left(1-\sqrt[3]{1-6\alpha}-2(1-6\alpha)^{-\frac{2}{3}}4\varepsilon\right)n-15\varepsilon n\\
&\geq \left(3-3\sqrt[3]{1-6\alpha}-\left(15+6(1-\beta_0)^{-2}\right)\varepsilon\right)n\\
&> \left(3-3\sqrt[3]{1-6\alpha}-39\varepsilon\right)n.
\end{align*}
where the first inequality follows from $\alpha -4\varepsilon<\frac{9}{2}\beta_0^3$ and $t\geq \frac{1}{\varepsilon}$, the second inequality follows from Lagrange's Mean Value Theorem, the third inequality follows from $\alpha \leq \frac{9}{2}\beta_0^3$, and the last inequality follows from $1-\beta_0>\frac{1}{2}$.

Thus, the theorem holds.
\qed

\begin{thm}\label{main3}
Let $\mathcal{L}_{n,k}^{(3)}$ be the family of all tight linear forests of order $n$ with $k$ edges in 3-graphs. For $0<c<1$ and $n$ sufficiently large,
\[
ex_3(n;\mathcal{L}_{n,cn}^{(3)})=\max\left\{\binom{cn}{3},
\binom{n}{3}-\binom{n-\frac{cn}{3}}{3}\right\}+o(n^3).
\]
\end{thm}
\pf Let
\begin{align*}
m=\max\left\{\frac{c^3n^3}{6},\frac{n^3}{6}-\frac{\left(1-\frac{c}{3}\right)^3n^3}{6}\right\}.
\end{align*}
Firstly, we prove that $ex_3(n;\mathcal{L}_{n,cn}^{(3)})\leq m+\gamma n^3$, where $\gamma>0$ is a small constant. Let $G$ be a graph with $m+\gamma n^3$ edges. If $c\geq 3\beta_0$, then $e(G) = \frac{c^3n^3}{6}+\gamma n^3 > \frac{9}{2}\beta_0^3 n^3$. By Theorem \ref{main-3}, we have
\[
l(H)\geq \left(\sqrt[3]{6\cdot\frac{c^3}{6}+\gamma}-o(1)\right)n \geq cn.
\]
 If $c< 3\beta_0$, then $e(G) = \frac{n^3}{6}-\frac{\left(1-\frac{c}{3}\right)^3n^3}{6}+\gamma n^3 < \frac{9}{2}\beta_0^3 n^3$. By Theorem \ref{main-3}, we have
\[
l(H)\geq 3\left(1-\sqrt[3]{1-6\left(\frac{1-\left(1-\frac{c}{3}\right)^3}{6}+\gamma\right)}-o(1)\right)n \geq cn.
\]
Thus, we concluded that
\[
ex_3(n;\mathcal{L}_{n,cn}^{(3)})\leq m+o(n^3) = \max\left\{\binom{cn}{3},
\binom{n}{3}-\binom{n-\frac{cn}{3}}{3}\right\}+o(n^3).
\]

On the other hand, the constructions $K_{cn}^{(3)}\cup E_{(1-c)n}^{(3)}$ and  $K_{(cn-1)/3}^{(3)}\vee E_{(1-c/3)n+1/3}^{(3)}$ show that
\[
ex_3(n;\mathcal{L}_{n,cn}^{(3)})\geq \max\left\{\binom{cn}{3},
\binom{n}{3}-\binom{n-\frac{cn}{3}}{3}\right\}+o(n^3).
\]

Thus, we complete the proof.
\qed

Combining Theorem \ref{mainr} and \ref{main3}, we prove the Theorem \ref{main2}.

\section{Concluding Remarks}

In this paper, we propose a conjecture for Tur\'{a}n number of a family of tight linear forests.
We verified it for parts of the dense hypergraph case.
As we have proved, on one hand, the proposed conjecture implies the Erd\H{o}s Matching Conjecture directly.
On the other hand, by weak  regularity lemma for hypergraphs, the Erd\H{o}s Matching Conjecture implies that the proposed
conjecture is asymptotically true in the dense case.

Since a tight linear forest on $n$ vertices with $n-r+1$ edges is exactly a Hamilton tight path.
We can view the tight linear forest as an intermedia concept between matching and Hamilton tight cycle. Thus, it seems that tight linear forests are natural generalizations of matchings.
It is an interesting problem to determine the exact Tur\'{a}n number of this family of  tight linear forests.

\noindent{\bf Acknowledgements.}
The work was supported by National Natural Science Foundation of China (No. 11701407, No. 11671296, No. 61502330).


\begin{thebibliography}{00}

\bibitem{bushaw}
N. Bushaw and N. Kettle, Tur\'{a}n numbers of multiple paths and equibipartite forests,
{\it Combin. Probab. Comput.} {\bf 20} (2011) 837--853.


\bibitem{bollobas2}
B. Bollob\'{a}s, D. E. Daykin, and P. Erd\H{o}s, Sets of independent edges of a hypergraph, {\it Quart. J. Math. Oxford Ser.}
(2) {\bf 21} (1976) 25--32.

\bibitem{chung}
F. R. K. Chung, Regularity lemmas for hypergraphs and quasi-randomness, {\it Random Structures Algorithms} (2) {\bf 2} (1991) 241--252.

\bibitem{erdos2}
P. Erd\H{o}s and T. Gallai, On maximal paths and circuits of graphs, {\it Acta Mathematica Academiae Scientiarum Hungarica} (3-4) {\bf 10} (1959) 337--356.

\bibitem{erdos3}
P. Erd\H{o}s, C. Ko, and R. Rado, Intersection theorems for systems of finite sets, {\it Quart. J. Math. Oxford Ser.} (2) {\bf 12} (1961) 313--320.

\bibitem{erdos4}
P. Erd\H{o}s, A problem on independent $r$-tuples, {\it Ann. Univ. Sci. Budapest. E\"{o}tv\"{o}s Sect. Math.} {\bf 8} (1965) 93--95.

\bibitem{frankl0}
P. Frankl, V. R\"{o}dl, and A. Ruci\'{n}ski, On the maximum number of edges in a triple system not containing
a disjoint family of a given size, {\it Combin. Probab. Comput.} {\bf 21} (2012) 141--148.

\bibitem{frankl1}
P. Frankl, Improved bounds for Erd\H{o}s' matching conjecture, {\it J. Combin. Theory Ser. A} (5) {\bf 120} (2013) 1068--1072.

\bibitem{frankl2}
P. Frankl, On the maximum number of edges in a hypergraph with given matching number, {\it Discrete Appl. Math.} {\bf 216} (2017) 562--581.

\bibitem{han10}
H. H\`{a}n and M. Schacht, Dirac-type results for loose Hamilton cycles in uniform hypergraphs,
{\it J. Combin. Theory Ser. B} (3) {\bf 100} (2010) 332--346.

\bibitem{huang}
H. Huang, P. S. Loh, and B. Sudakov, The size of a hypergraph and its matching number,
 {\it  Combin. Probab. Comput.} (3) {\bf 21} (2012) 442--450.

\bibitem{keevash}
P. Keevash, Hypergraph Tur\'{a}n problems, {\it Surveys in combinatorics} {\bf 392} (2011) 83--140.

\bibitem{lidicky}
B. Lidicky, H. Liu, and C. Palmer, On the Tur\'{a}n number of forests, {\it Electron. J. Combin.} (2) {\bf 20} (2013) \#P62.

\bibitem{ningbo}
B. Ning and J. Wang, The formula for Tur\'{a}n number of spanning linear forests,  arXiv:1812.01047, 2018.

\bibitem{szem}
E. Szemer\'{e}di, Regular partitions of graphs, Probl\`{e}mes combinatoires et th\'{e}orie des graphes, {\it Colloq Int CNRS} {\bf 260} (1978), 399--401.

\bibitem{sidorenko}
A. Sidorenko, What we know and what we do not know about Tur\'{a}n numbers,
{\it Graphs Combin.} {\bf 11} (1995) 179--199.

\bibitem{turan1}
P. Tur\'{a}n, On an extremal problem in graph theory (in Hungarian), {\it Mat. Fiz.
Lapok} {\bf 48} (1941) 436--452.

\bibitem{wj}
J. Wang and W. Yang, The Tur\'{a}n number for spanning linear forests,
{\it Discrete Appl. Math.} (2018) https://doi.org/10.1016/j.dam.2018.07.014.

\bibitem{zhang}
L. Yuan and X. Zhang, The Tur\'{a}n number of disjoint copies of paths,
{\it Discrete Math.} (2) {\bf 340} (2017) 132--139.
\end{thebibliography}
\end{document}